\documentclass[a4paper]{amsart}

\usepackage{amsfonts,amsthm}
\usepackage{amssymb,amsmath,graphics}
\usepackage{latexsym}
\usepackage{eucal}

\makeatletter
\@addtoreset{equation}{section}\makeatother

\newtheorem{theo}{Theorem}[section]
\newtheorem{lem}[theo]{Lemma}
\newtheorem{prop}[theo]{Proposition}

\def\IsomPlusH{PSL_2(\rr)}
\def\IsomPlusH{\mathop{\mathrm{Isom}}^+(\hh^2)}
\def\tr{\mathop{\mathrm{tr}}\nolimits}
\def\idroit{\mathrm{i}}
\def\e#1{\mathrm{e}^{#1}}
\def\eintheta{e^{\idroit n\theta}}

\newcommand{\mc}{\mathcal}
\newcommand{\rr}{\mathbb{R}}
\newcommand{\nn}{\mathbb{N}}
\newcommand{\cc}{\mathbb{C}}
\newcommand{\hh}{\mathbb{H}}
\newcommand{\zz}{\mathbb{Z}}

\newcommand{\la}{\lambda}
\newcommand{\eps}{\epsilon}

\newcommand{\pl}{\partial}
\newcommand{\x}{\times}

\newcommand{\cjd}{\rangle}
\newcommand{\cjg}{\langle}

\newcommand{\demi}{\frac{1}{2}}
\newcommand{\ndemi}{\frac{n}{2}}
\newcommand{\tra}{\textrm{Tr}}

\newcommand{\zerov}{\mathop{\hbox{\rm 0-vol}}\nolimits}

\newcommand{\Proof}{\noindent\textsl{Proof}: }

\title{The determinant of the Dirichlet-to-Neumann map for surfaces with boundary}
\author{Colin Guillarmou}
\address{Laboratoire J. Dieudonn\'e\\
CNRS, Universit\'e de Nice\\ Parc Valrose\\
06100 Nice\\ France}     
\email{cguillar@math.unice.fr}
\author{Laurent Guillop\'e}
\address{Laboratoire J. Leray\\
CNRS, Universit\'e de Nantes\\ 2, rue de la Houssini\`ere\\
BP 92208,  44322 Nantes cedex 03\\ France}     
\email{laurent.guillope@univ-nantes.fr}

\begin{document}

\begin{abstract}
For any orientable compact surface with boundary, we 
compute the regularized determinant of the Dirichlet-to-Neumann (DN) map
in terms of particular values of dynamical zeta functions
by using natural uniformizations, one due to Mazzeo-Taylor, the  other 
to Osgood-Phillips-Sarnak.
We also relate in any dimension the DN map for the Yamabe operator to 
the scattering operator for a conformally compact related problem
by using uniformization.
\end{abstract}

\maketitle

\section{Introduction}
Let $(\bar{X},\bar{g})$ be a connected compact Riemannian manifold with
boundary, then the Dirichlet-to-Neumann (DN) map is the map
\[\mc{N}: C^\infty(\pl\bar{X})\to C^\infty(\pl\bar{X})\] 
defined by the following problem: let $f\in C^{\infty}(\pl\bar{X})$ and let $u\in C^{\infty}(\bar{X})$ 
be the solution of 
\[ \Delta_{\bar{g}}u=0 , \quad u|_{\pl\bar{X}}=f,\]
then if $\pl_n$ is the interior pointing vector field which is normal to $\pl\bar{X}$, we set 
\[\mc{N}f:= -\pl_n u|_{\pl\bar{X}}.\]
It is well-known that $\mc{N}$ is an elliptic self-adjoint
pseudo-differential operator on $\pl\bar{X}$ with principal symbol
$|\xi|_{h_0},\xi\in T^*\pl\bar X$, if $h_0:=g|_{T\pl\bar{X}}$ (see
\cite[7.11]{DNmap} for example). It is then possible to define its
determinant by the Ray-Singer method \cite{RS}.  Indeed, if $A$ is an
elliptic self-adjoint pseudo-differential operator of order $p>0$ with
positive principal symbol, we can set, following \cite{Se,RS,KV},
\[\det(A)=e^{-\pl_s\zeta_A(0)}, \quad \zeta_A(s)=\tra(A^{-s})\]
where $\zeta_A(s)$ is a priori defined for $\Re(s)\gg 0$ but has a
meromorphic extension to $\cc$ with no pole at $s=0$. If we apply this to
$A=\mc{N}$, we obtain $\det(\mc{N})=0$ since $\ker\mc{N}\not=0$, indeed
by Green's identity one has 
$\int_{\pl\bar{X}}\mc{N}f.f{\rm dvol}_{h_0}=\int_{\bar{X}}|\nabla u|^2{\rm dvol}_{\bar{g}}$
thus $\ker\mc{N}=\rr$ is the space of constant functions on $\pl\bar{X}$.
We then have to modify the definition of $\det(\mc{N})$: if $\Pi$ is the orthogonal
projection in $L^2(\pl\bar{X},{\mathrm dvol}_{h_0})$ onto the kernel
$\mathrm{ker} \mc{N}$, we take $\det'\mc{N}$ defined by
\[
\mbox{$\det'$}(\mc{N})=e^{-\pl_s\zeta^*_{\mc{N}}(0)}, \quad
\zeta_{\mc{N}}^*(s)=\tra(\mc{N}^{-s}(1-\Pi))
\]
which is well defined as before.
To compute $\det'(\mc{N})$, we first show the
\begin{theo}\label{confinv}
If $(\bar{X},\bar{g})$ is a Riemannian surface with boundary and if $\ell_{\bar{g}}(\pl\bar{X})$ is the length of the boundary $\pl\bar{X}$ for the metric $\bar{g}$, the value $\det'(\mc{N})/\ell_{\bar{g}}(\pl\bar{X})$
is a conformal invariant of the conformal manifold with boundary $(\bar{X},[\bar{g}])$.
\end{theo}
Note that this was proved in the case where $\bar{X}$ is a topological disc by Edward-Wu \cite{EW}.\\

Consequently, it is sufficient to study the case of particular conformal representative in the conformal class, 
that is to use uniformization. 

The first natural uniformization we will use 
has been proved by Mazzeo and Taylor \cite{MT}, it picks a complete constant negative curvature 
metric in the conformal class:  
indeed, they show that there exists a unique conformally compact metric $g$ on the
interior $X$ of $\bar{X}$ such that $g$ has curvature $-1$ and $g$ is
conformal to $\bar{g}$. The manifold $(X,g)$ is then isometric to an
infinite volume quotient $X\simeq \Gamma\backslash\hh^{2}$ of the
hyperbolic plane by a convex co-compact group of isometries.
We use this uniformization to compute $\det'(\mc{N})$, although the
DN map in this case does not really make sense, but
instead we have the scattering operator.

Before stating the result, we need to recall a few definitions about
Riemann surfaces and their Selberg (resp. Ruelle) zeta function. Let
$\Gamma\subset \IsomPlusH$ be a Fuchsian subgroup with only hyperbolic
elements (i.e. fixing $2$ points at the boundary of $\hh^2$), the quotient
$X=\Gamma\backslash\hh^2$ is a geometrically finite complete hyperbolic
manifold.  We recall that any $\gamma\in \Gamma$ is conjugated to the
dilation $z\to e^{\ell(\gamma)}z$, with translation length
$\ell(\gamma)\in\rr^+$ in the hyperbolic half-plane model $\hh^2=\{z\in\cc;
\Im(z)>0\}$, note that the set $[\Gamma]$ of primitive conjugacy classes of
$\Gamma$ is in one-to-one correspondence with the set $[C]$ of primitive
closed oriented geodesics $c$, the length of the closed geodesic $c$
corresponding to $\gamma$ being equal to $\ell(\gamma)$.  There is a dynamical
Ruelle type zeta function defined by the formula\footnote{The original zeta
  function of Ruelle was actually defined by the inverse of this one, we
  prefer to use the convention of Fried \cite{FR1}.}
\begin{equation}\label{ruelle}
R_\Gamma(\la):=\prod_{[\gamma]\in [\Gamma]}(1-e^{-\la \ell(\gamma)})\end{equation}
and the Selberg zeta function 
\begin{equation}\label{selberg}
Z_\Gamma(\la):=
\prod_{k\in\nn_0}R_\Gamma(\la+k).
\end{equation} 
These products converge for
$\Re(\la)>\delta$ where $\delta\in [0,1]$ is the exponent of convergence of
the Poincar\'e series of $\Gamma$, equal to $1$ only if $\Gamma$ is
cocompact. Moreover they admit an analytic extension\footnote{For the
  compact case, this is a consequence of Selberg trace formula, here this
  follows from Fried \cite{FR1} and Patterson-Perry \cite{PP} for
  instance.} to $\cc$ and verify the identity
\begin{equation}\label{RvsZ}
R_\Gamma(\la)=Z_\Gamma(\la)/Z_\Gamma(\la+1).
\end{equation}
 
We are able to compute the determinant of the 
DN map $\mc{N}$ using the uniformization of \cite{MT}:
\begin{theo}\label{detDN}
  Let $(\bar{X},\bar{g})$ be a smooth compact orientable connected
  Riemannian surface of Euler characteristic $\chi(\bar{X})$, with boundary
  $\pl\bar{X}$ of length $\ell(\pl\bar{X})$, and let $\mc{N}$ be the
  Dirichlet-to-Neumann operator of $\Delta_{\bar{g}}$ on $\pl\bar{X}$. 
  Let $g$ be the unique up to isometry, infinite volume, complete
  hyperbolic metric $g$ on $X$ conformal to $\bar{g}_{X}$ and let
  $\Gamma\subset \IsomPlusH$ the geometrically finite Fuchsian group such
  that $(X,g)$ is isometric to the space form $\Gamma\backslash\hh^{2}$.
  If we denote by $R_\Gamma(\la)$ the Ruelle zeta function of $\Gamma$, we
  have
\[\frac{\det'(\mc{N})}{\ell(\pl\bar{X})}=
\begin{cases}
1 & \hbox{if $\chi(\bar{X})=1$},\\
\ell(\gamma)/\pi & \hbox{if $\chi(\bar{X})=0$},\\
\left[(2\pi\la)^{\chi(\bar{X})-1}R_{\Gamma}(\la)\right]_{|\la=0}/\chi(\bar{X}) & \hbox{if $\chi(\bar{X})<0$}.
\end{cases}
\]
In the second case $\chi(\bar X)=0$, the group $\Gamma$ is
cyclic elementary, generated  by the hyperbolic isometry $\gamma$ with
translation length $\ell(\gamma)$, length of the unique closed geodesic of
the cylinder $ X\simeq\Gamma\backslash\hh^2$.
\end{theo}

The proof of this theorem is based on a functional equation for Selberg
zeta function for convex co-compact groups obtained in previous work
\cite{G} and the observation that the DN map for
$(\bar{X},\bar{g})$ is, modulo constant, the scattering operator $\mc{S}(\la)$ of the
uniformized non-compact manifold $\Gamma\backslash\hh^2$ at the parameter
value $\la=1$ : this is discussed in more generality at the end of the
introduction.  We emphasize that the Theorem holds even when the boundary has more than 
one connected component, an important fact that we need in the proof
being that $\ker \mc{N}$ is always equal to the space of constants 
and not the locally constant functions.
\medskip

\textsl{Remark $1$:} Note that for odd dimensional closed hyperbolic manifolds $X=\Gamma\backslash \hh^{d}$, 
the value $|R_{\Gamma}(0)|$ for some acyclic  representation 
of the $\pi_1$ of the unit tangent bundle $SX$
is the Reidemeister torsion of $SX$ by a result of Fried \cite{Fr}.\\

\textsl{Remark $2$:}
It is also worth to say that the proof shows that $0$ is
always a resonance of multiplicity $1$ with resonant state $1$ for the
Laplacian on any convex co-compact surface except when it's a cylinder
where it is then of multiplicity $2$; as a byproduct it also gives 
the exact order of vanishing of $R_{\Gamma}(\la)$ at $\la=0$, which was not apparently known in that case.\\

The next natural uniformizations for oriented compact surfaces with boundary 
are given by Osgood-Phillips-Sarnak \cite{OPS} (see also Brendle
\cite{BR}), they are of two types: each conformal class of a metric on an oriented compact surface with boundary 
has a unique  
\begin{itemize}
\item metric with constant Gauss curvature and with totally
  geodesic boundary. 
\item flat metric with constant geodesic curvature boundary.   
\end{itemize}  
The Gauss curvature $K$ on $\bar{X}$ and the geodesic curvature $k$ of the boundary $\pl\bar{X}$ 
are linked through the Gauss-Bonnet formula
\[\int_{\bar{X}}K {\rm dvol}_g +\int_{\pl\bar{X}} k {\rm d\ell}_g=2\pi \chi(\bar{X}).\]
The flat uniformization has been used by Edward-Wu \cite{EW} to show that $\det'(\mc{N})=\ell(\pl\bar{X})$ for a topological disc (\emph{i.e.} $\chi(\bar{X})=1$), their explicit computation is possible thanks to 
the circular symmetry of the uniformized flat disc with constant geodesic curvature. 
In a similar way, we give in the appendix the explicit computation for
the planar annulus whose boundary is the union of two
concentric circles and show that it fits with the value found in Theorem
\ref{detDN} for the hyperbolic cylinder conformal to this annulus. 
In the case $\chi(\bar X)<0$, the computation of $\det'(\mc{N})$ in terms of geometric quantities 
by using the flat uniformization does not seem apparent at all. 
As for the constant curvature with geodesic boundary uniformization, the topological disc (i.e. $\chi(\bar{X})>0$)
is uniformized by a half-sphere of curvature $+1$, the topological cylinder (i.e. $\chi(\bar{X})=0$) by a flat
cylinder $[0,L]\x S^1$ and in both cases, the value $\det'(\mc{N})$ can be easily computed using decomposition in
spherical harmonics of the Laplacian, essentially like for the flat uniformization.    
However, when $\chi(\bar{X})<0$, the constant curvature $-1$ uniformization with totally geodesic boundary 
appears to be more useful to compute $\det'(\mc{N})$. Indeed it yields a metric on 
$\bar{X}$ which is isometric to $G_0\backslash\hh^2$ for some discrete group
  $G_0$ of isometries of $\hh^2$ (containing symmetries of order $2$) and the double of $G_0\backslash \hh^2$ 
  along the boundary is the closed hyperbolic surface $M:=G\backslash \hh^2$ where $G:=
G_0\cap\IsomPlusH$ is the index $2$ subgroup of orientation preserving isometries of $G_0$.
The Mayer-Vietoris formula for determinants by Burghelea-Friedlander-Kappeler \cite{BFK} reads in this case
\[\frac{\det'(\mc{N})}{\ell(\pl\bar{X})}=-\frac{1}{2\pi\chi(\bar{X})}\frac{\det'(\Delta_{G\backslash\hh^2})}{(\det(\Delta_{G_0\backslash\hh^2}))^2}\]
where $\Delta_{G_0\backslash \hh^2}$ is the Dirichlet realization of the Laplacian on $G_0\backslash \hh^2$.
We are thus interested in the value of the regularized determinants of these Laplacians. The determinant $\det'(\Delta_{G\backslash\hh^2})$ has been computed by Sarnak and Voros \cite{Sa,V} in terms of the derivative at $\la=1$ of the Selberg zeta function $Z_G(\la)$ defined by (\ref{selberg}) with $\Gamma=G$. Using a trace formula of \cite{Gu}, 
we prove a similar formula for $\det(\Delta_{G_0\backslash\hh^2})$ in terms of a Selberg zeta function 
$Z_{G_0}(\la)$ at $\la=1$, where the natural Ruelle and Selberg zeta functions $R_{G_0}(\la),Z_{G_0}(\la)$ for this case
with boundary are defined as follows (see \cite[Section 5]{Gu}): let
$\ell_{1},\dots,\ell_{N}$ be the lengths of the geodesic boundary
components of $\bar{X}$ and let $[C]$ be the set of primitive oriented closed
geodesics $c$ of length $\ell_c$ and with $n_c$ geometric reflections (according
to the geometric optic law) on $\pl\bar{X}$, then the zeta
functions\footnote{The function we use is actually the square root of that
  of \cite{Gu}.} are defined by the following products:
\begin{equation}
\label{SelbergG0}
\begin{split}
R_{\partial\bar X}(\la):=\prod_{j=1}^N(1-e^{-\la\ell_{j}})^2,&\quad
R_{G_0}(\la):=\prod_{c\in[C]}
(1-(-1)^{n_c}e^{-\la\ell_c})(1-e^{-(\la+1)\ell_c}),\\
Z_{G_0}(\la)&:=\prod_{k\in\nn_0}
R_{\partial\bar X}(\la+2k)
R_{G_0}(\la+2k).
\end{split}
\end{equation}
We show the 
\begin{theo}\label{detDirichlet}
  Let $(\bar{X},\bar{g})$ be a compact oriented surface with boundary, with negative Euler 
  characteristic $\chi(\bar{X})$. Let $\bar{g}_0$ be the unique, up to isometry, constant negative curvature metric with totally geodesic     boundary on $\bar{X}$ and let $G_0\subset {\rm Isom}(\hh^2)$ be the discrete group such that 
  $(\bar{X},\bar{g}_0)$ is isometric to $G_0\backslash \hh^2$.
  Let $G=G_0\cap \IsomPlusH$ be the subgroup of
  $G_0$ of orientation preserving isometries and $Z_{G_0}(\la),Z_{G}(\la)$ be the associated Selberg zeta
  function of $G_0$ and $G$, then
\[\frac{\det'(\mc{N})}{\ell(\pl\bar{X})}=-\frac{Z'_G(1)}{(Z_{G_0}(1))^2}\frac{e^{\ell(\pl\bar{X})/4}}{2\pi\chi(\bar{X})}\]
Thus, if $\Gamma\backslash\hh^2$ is the uniformization of $(\bar{X},\bar{g})$ given by Theorem \ref{detDN}, then
\[
\left[\la^{\chi(\bar{X})-1}R_{\Gamma}(\la)\right]_{|\la=0}=-\frac{Z'_G(1)}{Z_{G_0}(1)^2}e^{\frac{
      \ell(\pl\bar{X})}{4}}(2\pi)^{-\chi(\bar{X})}.\]
\end{theo}

Although the products defining $R_{\Gamma}(\la),Z_{G}(\la),Z_{G_0}(\la)$ do not converge,
we can view the last identity of Theorem \ref{detDirichlet}
as a relation between length spectrum of $\Gamma\backslash\hh^2$ and $G_0\backslash\hh^2$, which does not appear obvious
at all. 
Let us also remark that the
determinant of the Laplacian on an hyperbolic compact surface
has different expressions with Selberg zeta values,
the Sarnak-Voros \cite{Sa,V} one related to the Fuchsian uniformization and
the McIntyre-Takhtajan \cite{MaTa} related to the Schottky uniformization.\\

In the last section we discuss in more generality (in higher dimension) the
relation between Dirichlet-to-Neumann map and scattering operator.  An
$(n+1)$-dimensional \emph{asymptotically hyperbolic manifold} $(X,g)$ is a
complete Riemannian non-compact manifold, which is the interior of a smooth
compact manifold with boundary $\bar{X}$ such that for any boundary
defining function $x$ of $\pl\bar{X}$ (\emph{i.e.} $\pl\bar{X}=\{x=0\}$ and
$dx|_{\pl\bar{X}}\not=0$), then $\bar{g}=x^2g$ is a smooth metric on
$\bar{X}$ such that $|dx|_{x^2g}=1$ on $\pl\bar{X}$.  The metric
$h_0:=\bar{g}|_{T\pl\bar{X}}$ induced on $\pl\bar{X}$ depends on $x$ and
another choice of $x$ yields a metric on $\pl\bar{X}$ conformal to $h_0$,
we thus define the \emph{conformal infinity} of $(X,g)$ as the conformal
class of $[h_0]$ on $\pl\bar{X}$.  There is a natural meromorphic family of
operators (defined in Section 2) $\mc{S}(\la)$ (for $\la\in\cc$) called
\emph{scattering operator}, acting on $C^{\infty}(\pl\bar{X})$, these are
elliptic conformally covariant pseudo-differential operators of order
$2\la-n$ with principal symbol $|\xi|^{2\la-n}_{h_0}$ where
$h_0=\bar{g}|_{T\pl\bar{X}}$ is a conformal representative of the conformal
infinity of $(X,g)$.  When $g$ is
Einstein, Graham and Zworski \cite{GZ} showed that $\mc{S}(n/2+k)$ for $k\in\nn$ are conformal
powers of the Laplacian on the boundary $\pl\bar{X}$, initially defined in \cite{GJMS}. 
Since $\mc{S}(\la)$ has order $1$ when $\la=(n+1)/2$ and the same principal symbol than a
Dirichlet-to-Neumann  map on the compact manifold $(\bar{X},\bar{g})$,
we may expect that it is realized as a DN map for an elliptic compact
problem with boundary.  We observe that when $g$ has constant scalar
curvature (for instance if $g$ is Einstein), then
$\mc{S}(\frac{n+1}{2})$ is the
Dirichlet-to-Neumann map of the conformal Laplacian on a whole class of smooth metric
$\bar{g}$ on $\bar{X}$, conformal to $g$, with $\bar{g}|_{T\pl\bar{X}}=h_0$
and with minimal boundary $\pl\bar{X}$. Conversely it is clear that there
is no constant curvature uniformization when $n+1>2$, but instead there is
a solution of a singular Yamabe problem, that is, for a given
$(\bar{X},\bar{g})$, there exists an asymptotically hyperbolic metric with constant
scalar curvature on the interior $X$ in the conformal class of $\bar{g}$.  The
existence and regularity of such a solution of this singular Yamabe problem
is due to Aviles-Mac Owen \cite{AMC}, Mazzeo \cite{Ma} and
Andersson-Chru\'sciel-Friedrich \cite{ACF}.  If $K$ is the mean curvature of $\pl\bar{X}$ for $\bar{g}$ and $\mc{N}$ is 
the DN map for the conformal Laplacian $P=\Delta_{\bar{g}}+{\mathrm Scal}_{\bar{g}}(n-1)/(4n)$, we show that $\mc{N}+(n-1)K/2$ is the value $\mc{S}((n+1)/2)$  for a complete manifold with constant negative scalar curvature,
conformal to $\bar{g}$ on the interior $X$ of $\bar{X}$. Note that 
$\mc{N}+(n-1)K/2$ is known to be the natural conformally covariant operator on the boundary 
associated to $P$, see \cite{Ch}.
 \medskip

\noindent\textbf{Acknowledgement}. We thank E. Aubry, P. Delano\"e,
M. Harmer, A. Hassell, R. Mazzeo and P. Perry for useful discussions and
for pointing out the right references.  C.G.  acknowledges support of NSF
grant DMS0500788, and french ANR grants JC05-52556 and JC0546063.

\section{Computation of $\det'(\mc{N})$ using Mazzeo-Taylor uniformization}

We now recall the definition of the scattering operator $\mc{S}(\la)$ on an asymptotically hyperbolic manifold 
$(X,g)$ of dimension $n+1$. From Graham-Lee \cite{GRL}, for any choice $h_0$ in the conformal infinity $[h_0]$,
such a metric can be written uniquely in a collar neighbourhood $[0,\eps)_x\x\pl\bar{X}$ 
of the boundary under the form
\begin{equation}\label{model}
g=\frac{dx^2+h(x)}{x^2}, \quad h(0)=h_0
\end{equation} 
for some smooth $1$-parameter family of metric $h(x)$ on $\pl\bar{X}$ ($x$ is a boundary defining function of $\pl\bar{X}$). If $h(x)$ has a Taylor expansion at $x=0$
with only even powers of $x$, then $g$ is called \emph{even} (see \cite{G1}). 
If $g$ is even and $f\in C^{\infty}(\pl\bar{X})$, $\Re(\la)\geq n/2$ and $\la\notin n/2+\nn$, then
$\mc{S}(\la)f:=c(\la)u_\la^+|_{\pl\bar{X}}\in C^{\infty}(\pl\bar{X})$ where\footnote{We changed the convention since in the literature, $u_{\la}^+|_{\pl\bar{X}}$ would be the scattering operator acting on $f$.} $c(\la)$ is the normalisation constant $c(\la):=2^{2\la-n}\Gamma(\la-\ndemi)/\Gamma(\ndemi-\la)$ and $u_\la,u_\la^\pm$ are defined 
by solving the Poisson problem \cite{GZ}
\begin{equation}\label{poissonpb}
(\Delta_g-\la(n-\la))u_\la=0, \quad u_\la= x^{n-\la}u^-_\la+ x^{\la}u_\la^+
, \quad u_\la^\pm\in C^{\infty}(\bar{X}),\quad u_\la^-|_{\pl\bar{X}}=f.\end{equation}

We see that $\mc{S}(\la)f$ depends on $g$ and on the choice of $x$ or equivalently on the choice of conformal representative $h_0=x^2g|_{T\pl\bar{X}}$ of the conformal infinity of $(X,g)$. 
Changing $h_0$ into $\hat{h}_0=e^{2\omega_0}h_0$ with $\omega_0\in C^\infty(\pl\bar{X})$ induces the scattering operator
\begin{equation}\label{confcov}
\hat{\mc{S}}(\la)=e^{-\la \omega_0}\mc{S}(\la)e^{(n-\la)\omega_0}.\end{equation}
From \cite{JSB,GZ}, $\mc{S}(\la)$ is holomorphic in the half plane $\{\la\in\cc;\Re(\la)>n/2\}$, 
moreover it is a pseudodifferential operator of order $2\la-n$
with principal symbol $|\xi|^{2\la-n}_{h_0}$ (thus elliptic) and it is 
self-adjoint when $\la \in(n/2,+\infty)$, which
makes its zeta regularized determinant well defined by \cite{KV}.
If the dimension $n+1$ is even, one shows easily that if $\hat{h}_0$ is 
conformal to $h_0$, the conformal relation \eqref{confcov} between the associated operators $\mc{S}(\la)$ and $\hat{\mc{S}}(\la)$
implies that $\det(\mc{S}(\la))=\det(\hat{\mc{S}}(\la))$, see \cite[Sec. 4]{G} for
instance.
\medskip

We are back to our case of surfaces (here $n=1$), thus 
let $(\bar{X},\bar{g})$ be a smooth Riemannian surface with boundary.
We first relate $\mc{N}$ to the scattering operator of an associated non-compact 
hyperbolic surface.
Let $\rho$ be a function that defines $\pl\bar{X}$ and such that $\bar{g}=d\rho^2+h_0+O(\rho)$ 
for some metric $h_0$ on $\pl\bar{X}$, so the normal vector field to the boundary is $\pl_n=\pl_\rho$ on $\pl\bar{X}$. Let
$g = \hat{\rho}^{-2} \bar{g}$ be the unique complete hyperbolic metric on the interior $X$ of $\bar{X}$, obtained by Mazzeo-Taylor \cite{MT}, where 
$\hat{\rho}=\rho +O(\rho^2)$ is some smooth function on $\bar{X}$, then $(X,g)$ is an asymptotically hyperbolic
manifold in the sense stated in the introduction.
Then $g$ is even since the metric
outside some compact is the metric on a hyperbolic funnel, that is $dr^2+\cosh^2(r)dt^2$
on $(0,\infty)_r\x (\rr/a\zz)_t$ for some $a>0$ (it suffices to set $x=e^{-r}$ to have a model form \eqref{model}). 
Therefore the geodesic function $x$ such that $g$ is like \eqref{model} 
implies $h(x)=h_0+O(x^2)$ and $x=\rho+O(\rho^2)$.
By studying the Poisson problem at energy $\la$ close to $1$ for $\Delta_g$, for any $f\in C^{\infty}(\pl\bar{X})$, there exists a unique $u_\la\in C^{\infty}(X)$ such that 
\begin{equation}\label{poisson}
 (\Delta_g -\la(1-\la))u_\la =0 , \quad u_\la\sim_{x\to 0} x^{1-\la}\Big(f+\sum_{j=1}^\infty x^{2j}f^-_{2j}(\la)\Big) 
+ x^{\la}\Big(c(\la)\mc{S}(\la)f +\sum_{j=1}^\infty x^{2j}f^+_{2j}(\la)\Big)\end{equation}
for some $f^\pm_{2j}(\la)\in C^{\infty}(\pl\bar{X})$ (we used evenness of the metric so that odd powers of $x$ are zeros, see \cite{GZ}).
In particular at $\la=1$ we have $u:=u_1\in C^{\infty}(\bar{X})$ and $\Delta_g u=0$ but 
$\Delta_g=\hat{\rho}^2\Delta_{\bar{g}}$ thus 
\[ \Delta_{\bar{g}} u=0, \quad u\in C^{\infty}(\bar{X}), \quad u=f- x\mc{S}(1)f +O(\rho^2),\]
but since $\pl_x=\pl_{\hat{\rho}}=\pl_{\rho}=\pl_n$ on $\pl\bar{X}$ we automatically get
\begin{lem}\label{dnvsscattering}
The Dirichlet-to-Neumann map $\mc{N}$ for $\Delta_{\bar{g}}$ is given by the scattering operator $\mc{S}(1)$ at  energy $1$
for the Laplacian $\Delta_g$ on the asymptotically hyperbolic surface $(X,g)$ conformal to $\bar{g}$,
where $\mc{S}(\la)$ is defined using the boundary defining function associated to the representative 
$h_0=\bar{g}|_{T\pl\bar{X}}$ of the conformal infinity $[h_0]$ of $(X,g)$.
\end{lem}

Taking a conformal metric $\bar{g}_1=e^{2\omega}\bar{g}$ on $\bar{X}$ gives a
Laplacian $\Delta_{\bar{g}_1}=e^{-2\omega}\Delta_{\bar{g}}$ and the normal
vector field to the boundary becomes $\pl_n=e^{-\omega_0}\pl_x$ where
$\omega_0=\omega|_{\pl\bar{X}}$. We deduce that the associated
Dirichlet-to-Neumann map $\mc{N}_1$ satisfies
$\mc{N}_1=e^{-\omega_0}\mc{N}$. 
\begin{theo}\label{changefordet}
Let $\bar{g}_0$ and $\bar{g}_1=e^{2\omega}\bar{g}_0$ be two conformally related metrics on a surface with boundary 
$\bar{X}$, and let $\mc{N}_0,\mc{N}_1$ be the respective Dirichlet-to-Neumann operators.  
Then $\det'(\mc{N}_0)/\ell_{\bar{g}_0}(\pl\bar{X})=\det'(\mc{N}_1)/\ell_{\bar{g}_1}(\pl\bar{X})$ where
$\ell_{\bar{g}_i}(\pl\bar{X})$ is the length of the boundary for the metric $\bar{g}_i$, $i=0,1$.
\end{theo}
\Proof   By the main formula of Paycha-Scott
\cite{PS} (see also \cite{OK}),
\begin{equation}
  \label{eq:PaSc}
{\det} '(\mc{N}_i)=\exp\Big({\rm TR}(\log(\mc{N}_i)(1-\Pi_i))\Big), \quad i=0,1
\end{equation}
where ${\rm TR}$ is the Kontsevich-Vishik canonical trace defined in \cite{KV}, $\log(\mc{N}_i)$ is defined 
by a contour integral (see \cite{KV,PS} for details), and $\Pi_i$ the orthogonal projection onto $\ker\mc{N}_i$ with respect
to the volume density on $\pl\bar{X}$ induced by $h_i:=\bar{g}_i|_{\pl\bar{X}}$, i.e.
the projection onto the constants for the volume density ${\rm dvol}_{h_i}$.
It is important to note that this formula holds (\emph{i.e.}
Guillemin-Wodzicki residue trace does not show-up in the formula) since the DN maps $\mc{N}_i$  
have regular parity in the sense of \cite[Sect. 2]{G} and thus $\log \mc{N}_i$ as well: indeed, take 
the Mazzeo-Taylor uniformization $g$ of $\bar{g}_0$ (which is the same than that of $\bar{g}_1$),
then Proposition 3.6 of \cite{G} shows that the scattering operator 
$\mc{S}(\la)$ associated to $g$ has regular parity in the sense of \cite[Sect. 2]{G} since the hyperbolic 
metric $g$ is even; this implies by using Lemma \ref{dnvsscattering} 
that $\mc{N}_i$ has regular parity for $i=1,2$.
If $g_t=e^{2t\omega} g_0$ is a conformal change and $h_t=g_t|_{\pl\bar{X}}=e^{2t\omega_0}h_0$ where $\omega_0=\omega|_{\pl\bar{X}}$, 
then the DN map for the metric $g_t$ is
unitarily equivalent to the self-adjoint operator $\mc{N}_t:=e^{-t\frac{\omega}{2}}\mc{N}_0e^{-t\frac{\omega}{2}}$ on $L^2(\pl\bar{X},{\rm dvol}_{h_0})$, with $L^2$ kernel projector 
\begin{equation}\label{pit}
\Pi_t=(\ell_{h_t}(\pl\bar{X}))^{-1}e^{t\frac{\omega}{2}}\otimes e^{t\frac{\omega}{2}} 
\end{equation}
where $\ell_{h_t}(\pl\bar{X})$ is the length of $\pl\bar{X}$ for the metric $h_t$.
First we have that 
\[\pl_t(\log(\mc{N}_t)(1-\Pi_t))=(\pl_t \mc{N}_t)\mc{N}_t^{-1}(1-\Pi_t)-\log(\mc{N}_t)\pl_t\Pi_t\]
where $\mc{N}_t^{-1}$ is the unique operator (modulo $\Pi_t$) satisfying $\mc{N}_t\mc{N}_t^{-1}=
\mc{N}_t^{-1}\mc{N}_t=1-\Pi_t$, that is
 $\mc{N}_t^{-1}=(1-\Pi_t)e^{t\frac{\omega_0}{2}}\mc{N}_0^{-1}e^{t\frac{\omega_0}{2}}(1-\Pi_t)$:
  indeed, multiplying on the left by $\mc{N}_t$ gives
\begin{eqnarray*}
\mc{N}_t(1-\Pi_t)e^{t\frac{\omega_0}{2}}\mc{N}_0^{-1}e^{t\frac{\omega_0}{2}}(1-\Pi_t)&=&e^{-t\frac{\omega_0}{2}}(1-\Pi_0)e^{t\frac{\omega_0}{2}}(1-\Pi_t)\\
&=&(1-\ell_{h_t}(\pl\bar{X})e^{-t\omega_0}\Pi_t)(1-\Pi_t)=1-\Pi_t
\end{eqnarray*}
and the same holds by multiplying on the right by $\mc{N}_t$.
Thus taking the log derivative of $\det'(\mc{N}_t)$ with respect to $t$ gives 
(by the same arguments than \cite[Sec. 4]{G}) that
\[\pl_t\log({\det}'(\mc{N}_t))=-{\rm TR}(\omega_0(1-\Pi_t))-\tra((\log(\mc{N}_t)\pl_t\Pi_t)\]
where $\tra$ is the usual trace. Using (\ref{pit}), we compute 
$\pl_t\Pi_t=\demi(\omega_0\Pi_t+\Pi_t\omega_0)+\pl_t(\ell_{h_t}(\pl\bar{X})^{-1})\Pi_t$ but 
since $\log(\mc{N}_t)\Pi_t=\Pi_t\log(\mc{N}_t)=0$ and the trace is cyclic, we have 
\[\tra((\log(\mc{N}_t)\pl_t\Pi_t)=\tra (\log(\mc{N}_t)\omega_0\Pi_t)=\tra(\Pi_t\log(\mc{N}_t))=0.\] 
Now ${\rm TR}(\omega_0)=0$ since the Kontsevich-Vishik trace of a differential operator is $0$ in 
odd dimension (see \cite{KV}), but the trace of a smoothing operator is the integral on the diagonal 
of its Schwartz kernel, therefore
\[
{\rm TR}(\omega_0\Pi_t)=\tra(\omega_0\Pi_t)=\frac{\int_{\pl\bar{X}}\omega_0 e^{t\omega_0}{\rm dvol}_{h_0}}
{\int_{\pl\bar{X}}e^{t\omega_0}{\rm dvol}_{h_0}}=\pl_t
\log\Big(\ell_{h_t}(\pl\bar{X})\Big)
\]
Then integrating in $t\in[0,1]$ we get the right law for the determinant
\qed\\

We now prove Theorem \ref{detDN}.
\medskip

\textsl{Proof of Theorem \ref{detDN}}: Since $\mc{S}(1)=\mc{N}$, we have to compute $\det'(\mc{S}(1))$. 
It is clear that the kernel of $\mc{S}(1)$
is one dimensional, composed of the constants, since it is the case for $\mc{N}$. 
 According to the main formula of Paycha-Scott \cite{PS} we have for $\la>1/2$ 
\[\det (\mc{S}(\la))=\exp\Big({\rm TR} (\log \mc{S}(\la))\Big), \quad 
{\det}'(\mc{S}(1))=\exp\Big({\rm TR}(\log(\mc{S}(1))(1-\Pi))\Big)\] 
where $\Pi$ is the projection onto the constants. 
To compute $\det'(\mc{S}(1))$, we shall analyze $\det\mc{S}(\la)$ in the neighborhood of
$\la=1$.
 
From \cite[Th. 1.3]{G} and the fact that $\Delta_g$
has no $L^2$ zero-eigenvalue, $\det(\mc{S}(\la))$ has no pole in a neighbourhood of 
$\la=1$ and is holomorphic near $\la=1$ with a zero of order $\nu_1$ where the multiplicity $\nu_{\la_0}$ for $\la_0\in\cc$ is defined by 
\begin{equation}\label{nula0}
\nu_{\la_0}:=-\tra \Big({\rm Res}_{\la=\la_0}
(\pl_\la\mc{S}(\la)\mc{S}^{-1}(\la))\Big)
\end{equation}

Let us now compute $\nu_1$.  We consider the largest integer $k$ such that
there exists a holomorphic (in $\la$) family of functions $u_\la$ in $L^2(\pl\bar{X})$
with $u_1\in\ker\mc{S}(1)$, and such that $\mc{S}(\la)u_\la=O((\la-1)^k)$.
This maximum is achieved for some $u_\la$, is positive and is exactly
$\nu_1$ by Gohberg-Sigal theory (see \cite{GS} or \cite{G2}). Thus there
exists a family of functions $u_\la$ on $\partial \bar X$, holomorphic in
$\la$, with $u_1\in\ker \mc{S}(1)$ such that
$\mc{S}(\la)u_\la=(\la-1)^{\nu_1}\psi+O((\la-1)^{\nu_1+1})$ for some
function $\psi\not=0$. Then setting $u_\la=u_1+(\la-1)v+O((\la-1)^2)$ we get
the equation
\[\mc{S}(\la)u_\la=(\la-1)(\mc{S}(1)v+\mc{S}'(1)u_1)+O((\la-1)^2) \]
which we multiply with $u_1$, integrate and use self adjointness of $\mc{S}(1)$ with $\rr u_1=\ker\mc{S}(1)$
to deduce 
\begin{equation}\label{slau1}
\cjg \mc{S}(\la)u_\la, u_1\cjd=
(\la-1)\int_{\pl\bar{X}} u_1\mc{S}'(1)u_1\textrm{ dvol}_{h_0}+O((\la-1)^2).\end{equation}
Recall that $u_1$ is constant since in the kernel of $\mc{N}$, but following the notation of Fefferman-Graham \cite[Th. 4.3]{FG}, we have a kind of $Q$ curvature defined by\footnote{They
actually define $Q:=S'(1)1$ where $S(\la):=c(\la)\mc{S}(\la)$, so clearly
$Q=-\mc{S}'(1)1$ since $\mc{S}(1)1=0$ and $c(1)=-1$.}
$Q:=-\mc{S}'(1)1$   and they prove the identity\footnote{We emphasize that their proof
 is only based on Green's identity and evenness of the metric expansion at the boundary. In particular 
 it includes the case of hyperbolic surfaces.}
\[ \int_{\pl\bar{X}}Q{\rm dvol}_{h_0}=-\zerov(X)\]
where $\zerov(X)$ is the renormalized volume (also called $0$-volume) of $X$, \emph{i.e.} the constant $V$ in the expansion
\[ {\rm Vol}(x>\eps)= c_0\eps^{-1}+ V+O(\eps) , \quad \textrm{ as }\eps \to 0.\] 
It is however proved, from Gauss-Bonnet formula, by Guillop\'e-Zworski \cite{GZ2} and Epstein \cite[Appendix]{PP}
that 
\begin{equation}\label{gaussbonnet}
\zerov(X)= -2\pi \chi(\bar{X})=2\pi(2g+N-2)
\end{equation}
where $\chi(\bar{X})$ is the Euler characteristic of $\bar{X}$, 
$g$ is its genus and $N$ the number of boundary components.
It follows that the coefficient of $(\la-1)$ in (\ref{slau1}) does not vanish if $\chi(\bar{X})\not=0$
and then $\nu_1=1$. 

Recall (see \cite{GZ}) that $\mc{S}(\la)$ is self-adjoint for $\la$ real,
Fredholm, analytic in $\la$ near $\la=1$ and invertible in a small pointed
disc (of radius $\eps>0$) centered at $1$, moreover $\mc{S}(1)$ has $0$ as
isolated eigenvalue of multiplicity $1$, then one can use Kato perturbation
theory \cite[VII,3]{K} to deduce that for $\la$ near $\la=1$,
the spectrum of the operator $\mc{S}(\la)$ near $0$ is an isolated eigenvalue of multiplicity $1$,
that we denote $\alpha(\la)$; moreover it is holomorphic in $\la$ near $1$ and
there is a holomorphic $L^2$ normalized associated eigenvector $w_\la$.  We
have $w_\la=w_1+O(\la-1)$ where $w_1=\ell(\pl\bar{X})^{-1/2}\in\ker
\mc{S}(1)$, $\ell(\pl\bar{X})$ being the length of the curve $\pl\bar{X}$
for the metric $h_0$, and we get the equation
\begin{equation}\label{sula}
\mc{S}(\la)w_\la=\alpha(\la)w_\la,\quad
 w_\la = 
\ell(\pl\bar{X})^{-1/2}+(\la-1) v+O((\la-1)^2), \quad \alpha(\la)=(\la-1)\beta+O((\la-1)^2)
 \end{equation}
 for some $v\in C^\infty(\pl\bar{X})$ and $\beta\in\rr$. Taking a Taylor
 expansion of \eqref{sula} yields
\[\mc{S}(1)v+\mc{S}'(1)w_1=\beta w_1 \]
where we used the notation $'$ for $\pl_\la$. Pairing as before with $w_1$ and using that $\mc{S}(1)$ is self adjoint 
and previous arguments with $\int _{\pl\bar{X}}Q{\rm dvol}_{h_0}=2\pi\chi(\bar{X})$ gives  
\begin{equation}\label{valuealpha1} 
\beta=-\frac{2\pi\chi(\bar{X})}{\ell(\pl\bar{X})}.
\end{equation}

Now let $\Pi_\la$ be the orthogonal projection onto $\ker (\mc{S}(\la)-\alpha(\la))$, and
we define the function
\[h(\la):=\exp\Big({\rm TR}(\log(\mc{S}(\la))(1-\Pi_\la))\Big).\]
It is analytic near $\la=1$ and $\Pi_1=\Pi$ thus the limit of $h(\la)$ at $\la=1$ is $h(1)={\det}'(\mc{S}(1))$ by (\ref{eq:PaSc}), the value we search to compute.
For $\la\in\rr$ close to $1$ but $\la\not=1$, we have $\log(\mc{S}(\la))\Pi_\la=\log(\alpha(\la))\Pi_\la$, 
thus using the first identity in (\ref{eq:PaSc}) and the fact that ${\rm TR}$ is the usual trace on finite 
rank operators, we obtain 
\[h(\la)=\exp({\rm TR}(\log \mc{S}(\la))\exp(-\log \alpha(\la))=\det(\mc{S}(\la))/\alpha(\la).\]
But  $\alpha(\la)=(\la-1)\beta(1+O(\la-1))$ by \eqref{sula},
which proves that
\[
{\det}'(\mc{S}(1))=\lim_{\la\to 1}\frac{\det \mc{S}(\la)}{\beta (\la-1)}.
\]
In \cite{G}, we proved the functional equation 
\[
\det(\mc{S}(\la))=\frac{Z_\Gamma(1-\la)}{Z_\Gamma(\la)}\exp\Big(-2\pi\chi(\bar{X})\int_0^{\la-\demi}t\tan(\pi t)dt\Big).
\]
which, following Voros \cite[Eq 7.24, 7.25]{V}, can be written under the form 
\begin{gather}
\begin{gathered}\label{functeq} 
\det(\mc{S}(\la))
=\frac{Z_\Gamma(1-\la)}{Z_\Gamma(\la)}\Big(\frac{(2\pi)^{1-2\la}\Gamma(\la)G(\la)^2}{\Gamma(1-\la)G(1-\la)^2}
\Big)^{-\chi(\bar{X})}
\end{gathered}
\end{gather}
where $G$ is the Barnes function (see \cite[Appendix]{V}) which satisfies in particular 
$\Gamma(z)G(z)=G(z+1)$ and $G(1)=1$.
Writing $Z_\Gamma(1-\la)=R_\Gamma(1-\la)Z_\Gamma(2-\la)$ in
\eqref{functeq} and using that $Z_\Gamma(\la)$ is holomorphic at $\la=1$ implies that 
\begin{eqnarray*}
{\det}'(\mc{S}(1))&=& -\beta^{-1}\left[\lim_{\la\to 1}\frac{R_{\Gamma}(1-\la)}{(1-\la)^{1-\chi(\bar{X})}}
\right]\lim_{\la\to 1}\left[\frac{(1-\la)(2\pi)^{1-2\la}\Gamma(1-\la)\Gamma(\la)G(\la)^2}{G(2-\la)^2}
\right]^{-\chi(\bar{X})}\\
&=& (2\pi)^{\chi(\bar{X})-1}\frac{\ell(\pl\bar{X})}{\chi(\bar{X})}\lim_{\la\to 1}\frac{R_{\Gamma}(1-\la)}{(1-\la)^{1-\chi(\bar{X})}}.
\end{eqnarray*}
and we are done when $\chi(\bar{X})<0$. 

If $\chi(\bar{X})=1$, $\bar{X}$ is a topological disc and the uniformization that puts a complete hyperbolic
metric on $X$ is the usual hyperbolic disc. The proof \cite{G} of the formula \eqref{functeq} remains true by
setting $Z_\Gamma(\la):=1$ and we can proceed as before where now $\Pi_\la$ is the projection on the 
constants $\Pi_{\la}=w_1\cjg w_1,.\cjd$ if $w_1$ is like above. We finally obtain
\[{\det}'(\mc{S}(1))=-\frac{\ell(\pl\bar{X})}{2\pi}\lim_{\la\to 1}\frac{\det(\mc{S}(\la))}{\la-1}=
\ell(\pl\bar{X})\]
and we are done for this case. Notice that it matches with the result of Edward-Wu \cite{EW}.

The last case $\chi(\bar{X})=0$ corresponds to the cylinder, whose interior
$X$ is uniformized by the cyclic elementary group $\Gamma=\cjg\gamma\cjd$,
with a unique closed geodesic of length $\ell$, the translation length of
the generator $\gamma$.
In other words, $X$ is conformal to the hyperbolic cylinder $H_\ell=(\rr_r\x
(\rr/\zz)_t, g=dr^2+\cosh^2(r)\ell^2dt^2)$ where $\ell>0$ is the length of
the unique closed geodesic $\{r=0\}$.  By Theorem \ref{changefordet}, it suffices
to compute it for the conformal representative of the boundary at infinity
$(|\sinh(r)|^{-2}g)|_{T\pl\bar{X}}=\ell^2 dt^2$ and the result will be
given by multiplying by $\ell_{h_0}(\pl\bar{X})/2\ell$.  The scattering
matrix for the conformal representative $\ell^2dt^2$ is computed in
\cite{GZ}, it is decomposable on the Fourier modes, the solution of
\eqref{poisson} for data $f=1$ on the boundary is
\begin{equation*}
\begin{split}
u_\la(r,t)=& |\sinh r|^{\la-1}F\Big(\frac{1-\la}{2},1-\frac{\la}{2},\frac{3}{2}-\la;-\sinh^{-2}(r))\Big)\\
&+\Big(\frac{\Gamma(\la/2)}{\Gamma((1-\la)/2)}\Big)^2\frac{\Gamma(\demi-\la)}{\Gamma(\la-\demi)}
|\sinh r|^{-\la}F\Big(\frac{\la}{2},\frac{\la+1}{2},\la+\demi;-\sinh^{-2}(r)\Big)
\end{split}
\end{equation*}
where $F$ is the hypergeometric function but since here we chose $x=|\sinh r|^{-1}$ (to have 
the right conformal representative on $\pl\bar{X}$) and since $F(a,b,c;0)=1$, this gives easily $\mc{S}(\la)1$:
\[\mc{S}(\la)1=2^{2\la-1}\Big(\frac{\Gamma(\la/2)}{\Gamma((1-\la)/2)}\Big)^2=\frac{\pi}{2}(1-\la)^2+O((1-\la)^3), \quad \la\to 1.\]
We can thus do the same reasoning as above, but now $\nu_1=2$, $w_\la=w_1$ and
$\beta=\pi/2$, we finally get, for the conformal representative $\ell^2dt^2$,
\[{\det}'(\mc{S}(1))=\lim_{\la\to 1}\frac{\det \mc{S}(\la)}{\pi(\la-1)^2/2}=\frac{2}{\pi}\lim_{\la\to 1}\frac{R_{\Gamma}(1-\la)}{(1-\la)^2}=\frac{2\ell^2}{\pi}\]
where the Selberg zeta function for this special case in the functional equation \eqref{functeq}
is $Z_{\Gamma}(\la)=\prod_{k\in\nn_0}(1-e^{-(\la+k)\ell})^2$
(see Prop 3.3 of Patterson \cite{Pa}), thus
here $R_{\Gamma}(\la)=(1-e^{-\la \ell})^2$. 
This gives the proof.    
For completeness, we will give another explicit computation of this case in the Appendix using
the flat annulus conformal to the hyperbolic cylinder $H_\ell$.
\qed\\

We remark that in the proof above, the fact that $\nu_1=1$ when $\chi(\bar{X})<0$
shows that $\det\mc{S}(\la)$ has a zero of order exactly $1$ at $\la=1$ and 
from (\ref{functeq}) we deduce that $Z_{\Gamma}(\la)$ has a zero of order exactly $-\chi(\bar{X})+1$, 
as well as the Ruelle function $R_{\Gamma}(\la)$.
It also implies that 
$\nu_{0}=-1$ since $\nu_{\la_0}=-\nu_{1-\la_0}$ for any $\la_0\in\cc$ (see comments after equation 1.1 of \cite{PP}).
Then, using also that $\la=1$ is not a pole of the resolvent \footnote{The resolvent $R(\la)$ is an analytic family of operators acting on $L^2(X)$ if $\Re(\la)>1$, it admits a meromorphic continuation to $\la\in\cc$
 as an operator mapping $C_0^\infty(X)$ to $C^{\infty}(X)$ by a result of Guillop\'e-Zworski \cite{GZ1}.} $R(\la):=(\Delta-\la(1-\la))^{-1}$
since $0$ is not an $L^2$ eigenvalue of $\Delta_g$, we deduce
from Theorem 1.1 of \cite{G2} that $\la=0$ is a pole of order $1$, with 
residue of rank $1$, of the meromorphic extension of $R(\la)$ to $\cc$.
 
\section{Computation of $\det'(\mc{N})$ using Osgood-Phillips-Sarnak uniformization}
In \cite{OPS}, Osgood, Phillips and Sarnak (see also Brendle \cite{BR}) 
proved that in each conformal class of metrics on a compact surface with boundary,  
\begin{itemize}
\item there is a unique representative which is flat, with constant geodesic curvature boundary,
\item there is a unique representative which has constant curvature and totally geodesic boundary. 
\end{itemize}
The flat uniformization has been used by Edward-Wu \cite{EW} to compute $\det'(\mc{N})$ for a topological disc (i.e. the case $\chi(\bar{X})=1$),
they found the same result than in our Theorem \ref{detDN}. One can do the same for a topological cylinder, 
it is uniformized as a flat annulus $\{z\in\cc; 1\leq |z|\leq \rho\}$ for some $\rho>1$ 
and it is possible to compute $\det'(\mc{N})$, we do the calculation in the Appendix and show that it fits
with the value in Theorem \ref{detDN}. If $\chi(\bar{X})<0$, there does not seem to be apparent
way to express $\det'(\mc{N})$ in terms of geometric invariants.
Thus, we use the constant curvature uniformization with geodesic boundary, in this case one obtains  
a compact hyperbolic surface with geodesic boundary $(\bar{X},\bar{g})$.
The surface $\bar{X}$ is uniformized so that $\bar{X}$ is isometric to
$G_0\backslash\hh^2$ where $G_0$ is a group of isometries of $\hh^2$,
containing some symmetries.  Associated to $G_0$, there is a natural
Selberg type zeta function $Z_{G_0}(\la)$ \cite{Gu}, defined in
\eqref{SelbergG0}.  Now, let $M=\bar{X}\sqcup\bar{X}$ be the manifold
obtained by gluing two copies of $\bar{X}$ at the boundary $\pl\bar{X}$,
then $M$ has smooth structure of surface with no boundary such that the
natural involution is smooth. We can extend the hyperbolic metric $\bar{g}$ on $M$ by
symmetry and the new metric, called $g$, is smooth on $M$ since the
structure of the metric $\bar{g}$ in Fermi coordinates $(r,t)\in[0,\eps)\x
(\rr/\zz)$ near each connected component $C=\{r=0\}$ of $\pl\bar{X}$ is
\[\bar{g}=dr^2+\ell^2\cosh^2(t)dt^2\]  
for  $\ell>0$ the length of $C$.
The manifold $M$  is isometric to the quotient $G\backslash\hh^2$ of the hyperbolic plane by the co-compact Fuchsian
group $G=G_0\cap\IsomPlusH$, the subgroup of index $2$ of direct isometries of $G_0$, we will call $(M,g)$ \emph{the double of $(\bar{X},\bar{g})$}.
Let us denote by $V:=\pl\bar{X}$ the boundary of $\bar{X}$,
the manifold $M\setminus V$ can be compactified canonically so that it corresponds
to two connected components isometric to $\bar{X}$, we will consider this manifold 
and will denote it $\bar{X}^2$ by abuse of notation. We show
\begin{theo}\label{detDirichlet1}
Let $(\bar{X},\bar{g})$ be an oriented surface with boundary with Euler characteristic $\chi(\bar{X})<0$.  
Let $G_0$ be the discrete subgroup of ${\rm Isom}(\hh^2)$ such that $(\bar{X},\bar{g})$ is conformal
to the hyperbolic surface with geodesic boundary $G_0\backslash\hh^2$ and  
let $G:=G_0\cap\IsomPlusH$ be the index $2$ subgroup of orientation preserving elements of $G_0$, so that
$G\backslash\hh^2$ is the closed hyperbolic surface realized by doubling $G_0\backslash\hh^2$ along the boundary. 
Then  
\[\frac{\det'(\mc{N})}{\ell(\pl\bar{X})}=-\frac{Z'_G(1)}{(Z_{G_0}(1))^2}\frac{e^{\ell(\pl\bar{X})/4}}{2\pi\chi(\bar{X})}\]
where $Z_G(\la), Z_{G_0}(\la)$ are the Selberg zeta function associated respectively to the group $G,G_0$, and
defined respectively in \eqref{selberg} and \eqref{SelbergG0}.
\end{theo}
\Proof By conformal invariance of $\det'(\mc{N})/\ell(\pl\bar{X})$ and using the constant negative curvature with totally geodesic boundary uniformization, it suffices to assume that $\bar{X}=G_0\backslash\hh^2$ as above. Following the notation preceding the Theorem, we let $M=G\backslash \hh^2$ be the double of $\bar{X}$ where $G=G_0\cap{\rm Isom}^+(\hh^2)$. 
We denote by $\Delta_{\bar{X}^2}$ the Laplacian on $M$ with
Dirichlet condition on the geodesic boundary $\pl\bar{X}$, that is the
direct sum $\Delta_{\bar{X}}\oplus\Delta_{\bar{X}}$ on the two copies of
$\bar{X}$ in $M$ where $\Delta_{\bar{X}}$ is the Dirichlet realization of the Laplacian on $\bar{X}$, then its spectrum is clearly the same than
$\Delta_{\bar{X}}$ but with double the multiplicity and thus
\begin{equation}\label{delta2}
\det (\Delta_{\bar{X}^2})=(\det(\Delta_{\bar{X}}))^2.\end{equation}
From the proof of Theorem $B^*$ of 
Burghelea-Friedlander-Kappeler \cite{BFK}, we get that, if $\Delta_M$ is the Laplacian on $(M,g)$ 
\begin{equation}\label{bfk}
\frac{{\det}'(\Delta_M)}{\det(\Delta_{\bar{X}^2})}=\frac{{\rm vol}(M)}{2\ell(\pl\bar{X})}{\det}'(\mc{N}) 
\end{equation}
where $\mc{N}$ is the Dirichlet-to-Neumann map on $\pl\bar{X}$ defined in the Introduction for either 
copy of $\bar{X}$ in $M$, ${\rm vol}(M)$ is the volume of $M$ for the hyperbolic metric $g$ (\emph{i.e.}
$-4\pi\chi(\bar{X})$ by Gauss-Bonnet formula), $\ell(\pl\bar{X})$ is the
length of the geodesic boundary  $\pl\bar{X}$.
But from Sarnak \cite{Sa} (see also D'Hoker-Phong \cite{DP}, Voros \cite{V})
we have (recall $\chi(M)=2\chi(\bar{X})$)
\begin{equation}\label{eta}
{\det}'(\Delta_M)=Z'_G(1)e^{-2\eta\chi(\bar{X})}, \quad \eta=2\zeta'(-1)-\frac{1}{4}+\demi\log(2\pi)
\end{equation}
with $\zeta$ the Riemann zeta function. 

We now need to compute the determinant $\det(\Delta_{\bar{X}})$ of the 
Dirichlet Laplacian on a hyperbolic surface with geodesic boundary $\bar{X}$, 
using the Selberg function $Z_{G_0}(\la)$. This can be done by methods of Sarnak \cite{Sa} and
a trace formula by Guillop\'e \cite[Prop 3.1]{Gu}, we show the
\begin{prop}\label{guillope}
If $\bar{X}=G_0\backslash \hh^2$ is a hyperbolic surface with geodesic boundary $\pl\bar{X}$ of length $\ell(\pl\bar{X})$ and $\Delta_{\bar{X}}$ is 
the Dirichlet Laplacian on $\bar{X}$, then
\[\det(\Delta_{\bar{X}}-\la(1-\la))=Z_{G_0}(\la)\Big( e^{\eta-\frac{\ell(\pl\bar{X})}{8\chi(\bar{X})}(1-2\la)+\la(1-\la)}\frac{(2\pi)^{\la-1}}{G(\la)^2\Gamma(\la)}\Big)^{-\chi(\bar{X})},\]
where $\chi(\bar{X})$ is the Euler characteristic of $\bar{X}$, $G(\la)$ is Barnes' function, 
$Z_{G_0}(\la)$ is the Selberg zeta 
function of \eqref{SelbergG0} and $\eta$ is the constant defined in \eqref{eta}. 
\end{prop}
\Proof It suffices to apply the proof of Sarnak \cite{Sa} (done in the case with no boundary) to the 
trace formula obtained in Proposition 3.1 of \cite{Gu}
\begin{equation*}
\begin{split}
\frac{Z'_{G_0}(\la)}{(2\la-1)Z_{G_0}(\la)}=&\tra(R_{\bar{X}}(\la)-R_{\bar{X}}(\la_0))+
\frac{Z_{G_0}'(\la_0)}{(2\la_0-1)Z_{G_0}(\la_0)}\\&-\chi(\bar{X})
\Big(\frac{\Gamma'(\la)}{\Gamma(\la)}-\frac{\Gamma'(\la_0)}{\Gamma(\la_0)}\Big)-\frac{\ell(\pl\bar{X})}{4}
\left(\frac{1}{2\la-1}-\frac{1}{2\la_0-1}\right),\end{split}
\end{equation*}
here $R_{\bar{X}}(\la):=(\Delta_{\bar{X}}-\la(1-\la))^{-1}$. 
Indeed one deduces from this the formula
\begin{equation}\label{CD}
  \det(\Delta_{\bar{X}}-\la(1-\la))=Z_{G_0}(\la) e^{-\frac{\ell(\pl\bar{X})}{4}\la+C\la(1-\la)+D}\Big(\frac{(\Gamma_2(\la))^2(2\pi)^{\la}}{\Gamma(\la)}\Big)^{-\chi(\bar{X})}\end{equation}
for some constants $C,D$ and where the digamma function $\Gamma_2$ is the
inverse of the Barnes function $G$ used in \eqref{functeq}. To compute the
constants, we consider the asymptotics as $\la\to+\infty$ of this identity.
First we get an asymptotic for the left hand side through the use of the
heat kernel small time asymptotic as in \cite{Sa}:
\begin{equation*}
\begin{split}
\log(\det(\Delta_{\bar{X}}-\la(1-\la)))=&-a_1\la(\la-1)\log(\la(\la-1))+a_1\la(\la-1)+2\sqrt{\pi}a_2(\la-\demi)\\
&+a_3\log(\la(\la-1))+o(1)\\
 =&-a_1\la^2\log\la+a_1\la^2+2a_1\la\log\la+2\sqrt{\pi} a_2\la+2a_3\log\la\\
 &-\sqrt{\pi}a_2-\demi a_1+o(1)
\end{split}
\end{equation*}
where $a_1,a_2,a_3$ are the heat invariant obtained in \cite[Appendix]{OPS} 
\[\tra(e^{-t\Delta_{\bar{X}}})=t^{-1}(a_1+a_2t^{\demi}+a_3t)+o(1)= \Big(-\demi\chi(\bar{X})-\frac{\ell(\pl\bar{X})}{8\sqrt{\pi}}t^{\demi}+
\frac{\chi(\bar{X})}{6}t\Big)+o(1)\] as $t\to 0$.
On the other hand, using Stirling formula and $\lim_{\la\to +\infty}Z_{G_0}(\la)=1$ in the right hand side of \eqref{CD}, we get
 \begin{equation*}
\begin{split}
2\log(\det(\Delta_{\bar{X}}-\la(1-\la)))=&-\chi(\bar{X})\Big(-(\la-1)^2\log(\la-1)+\frac{3}{2}(\la-1)^2-
\la\log\la+\la\\
&+\frac{2}{3}\log\la+\log(2\pi)-2\zeta'(-1)\Big)
 +C\la(1-\la)+D-\frac{\ell(\pl\bar{X})}{4}\la+o(1)\\
=&-\chi(\bar{X})\Big(-\la^2\log\la+\frac{3}{2}\la^2+\la\log\la-\la-\frac{1}{3}\log\la+\demi\log(2\pi)\\
&-2\zeta'(-1)\Big)-C\la^2+\Big(C-\frac{\ell(\pl\bar{X})}{4}\Big)\la+D+o(1).
\end{split}
\end{equation*}
Identifying the coefficient gives 
\[C=-\chi(\bar{X}), \quad D=\chi(\bar{X})\left(\demi\log (2\pi)-2\zeta'(-1)+\frac{1}{4}\right)+\frac{\ell(\pl\bar{X})}{8}\]
and this gives the desired formula.
\qed\\

The Theorem \ref{detDirichlet1} is now proved by combining \eqref{bfk}, \eqref{eta}, \eqref{delta2} and Proposition \ref{guillope} with 
$\la=1$ (note that $Z_{G_0}(\la)$ is holomorphic with no zero near $\la=1$, 
as $\Delta_{\bar X}$ has no zero eigenvalue).
\qed\\

Observe that $\det(\Delta_{\bar{X}}-\la(1-\la))^2=\det(\Delta_{\bar{X}^2}-\la(1-\la))$
where $\bar{X}^2$ is two disjoint copies of $\bar{X}$ and $\Delta_{\bar{X}^2}=\Delta_{\bar{X}}\oplus\Delta_{\bar{X}}$ 
the Dirichlet Laplacian. Similarly one has $\chi(M)=2\chi(\bar{X})$ thus one can use Sarnak's formula \cite{Sa}
to get
\[\frac{\det(\Delta_{M}-\la(1-\la))}{\det(\Delta_{\bar{X}^2}-\la(1-\la))}=
\frac{Z_{G}(\la)}{Z_{G_0}(\la)^2}e^{-\frac{ \ell(\pl\bar{X})}{4}(1-2\la)}.
\]
Moreover, if $\Gamma\backslash \hh^2$ is the Mazzeo-Taylor uniformization of $\bar{X}$, then Proposition \ref{guillope} and Theorem \ref{detDirichlet1} together give the formula  
\[\left[\la^{\chi(\bar{X})-1}R_{\Gamma}(\la)\right]_{|\la=0}=-\frac{Z'_G(1)}{Z_{G_0}(1)^2}e^{\frac{ \ell(\pl\bar{X})}{4}}(2\pi)^{-\chi(\bar{X})}.\]
This formula relates in a regularized way the  the length spectrum of $\bar{X}$ and 
that of the non-compact uniformization $\Gamma\subset \hh^2$.
It does not appear obvious to us how to obtain any relation between these spectrum by 
other methods.

\section{Dirichlet-to-Neumann for Yamabe operator}

We discuss now what is the higher dimensional version of Lemma \ref{dnvsscattering}.\\

First, let $(X,g)$ an asymptotically hyperbolic manifold with constant scalar curvature 
equal to ${\rm Scal}_g=-n(n+1)$. From Graham-Lee \cite{GRL}, there is boundary defining function
$x$ such that the metric near the boundary is $(dx^2+h(x))/x^2$ for some 1-parameter family of metric 
$h(x)$ on $\pl\bar{X}$. A straightforward computation gives  
\[{\rm Scal}_{g}=-n(n+1)= -n(n+1)+n x\pl_x\log(\det h(x))+x^2{\rm Scal}_{\bar{g}}\]
which implies that $\tr_{h_0}(h_1)=0$ if $h_0=h(0)$ and $h_1=\pl_xh(0)$.
The Poisson problem (\ref{poissonpb}) with initial data $f\in C^{\infty}(\pl\bar{X})$ can then be solved at $\la=(n+1)/2$ by results of \cite{GZ,G1} (for such $\la$, we do not need full evenness of $g$ but only $\tr_{h_0}(h_1)$). Since 
\[\Big(\Delta_{g}-\frac{n+1}{2}.\frac{n-1}{2}\Big)\hat{\rho}^{\frac{n-1}{2}}u=0, \quad
u\in C^{\infty}(\bar{X}), \quad u|_{\pl\bar{X}}=f 
\] 
is equivalent to solve the elliptic Dirichlet problem  
\begin{equation}\label{yamabe}
\Big(\Delta_{\bar{g}} +{\rm Scal}_{\bar{g}}\frac{n-1}{4n}\Big)u=0, \quad u|_{\pl\bar{X}}=f
\end{equation}
for $\bar{g}=x^2g$ by conformal covariance of the Yamabe operator, and since $\pl_n=\pl_x$ at the boundary, we deduce that $\mc{S}((n+1)/2)=\mc{N}$ where $\mc{N}$ is the DN map for the conformal Laplacian of $(\bar{X},\bar{g})$ and $\mc{S}(\la)$ is the scattering operator for $(X,g)$ with boundary defining function. Remark that $\pl\bar{X}$ is a minimal hypersurface of $(\bar{X},\bar{g})$ since $h_1=0$.\\

Conversely, let $(\bar{X},\bar{g})$ be an $(n+1)$-dimensional smooth compact Riemannian manifold with boundary,
then it is proved by Aviles-Mac Owen \cite{AMC} that there exists a complete metric
$g_0$ conformal to $\bar{g}$ on the interior $X$ and with negative constant scalar curvature ${\rm Scal}_{g_0}=-n(n+1)$. 
Moreover it is proved by Andersson-Chru\'sciel-Friedrich \cite[Th. 1.3]{ACF} (see also 
Mazzeo \cite{Ma}) that $g_0$ is asymptotically hyperbolic with $\log$ terms in the expansion,
more precisely let $\rho$ be a geodesic boundary defining function of $\pl\bar{X}$ for $\bar{g}$,
\emph{i.e.} $\bar{g}=d\rho^2+\bar{h}(\rho)$ for some $1$ parameter family of metric $\bar{h}(\rho)$ on $\pl\bar{X}$, we have 
\[ g_0=\frac{\bar{g}(1+\rho v+\rho^nw) }{\rho^2}=\frac{\bar{g}}{\hat{\rho}^2},\] 
with $v\in C^{\infty}(\bar{X}), w\in C^{\infty}(X)$ and $w$ having a polyhomogenous expansion
\[w(\rho,y)\sim \sum_{i=0}^\infty \sum_{j=0}^{N_i}u_{ij}\rho^i(\log \rho)^j\]
near the boundary, $N_i\in\nn_0$ and $u_{ij}\in C^{\infty}(\pl\bar{X})$.
Note that by Graham-Lee Lemma \cite{GRL}, there exists for $h_0:=\bar{h}(0)=\bar{g}|_{T\pl\bar{X}}$ a boundary defining function $x=\rho+O(\rho^2)$ such that
$g_0=(dx^2+h(x))/x^2$ near $\pl\bar{X}$ with $h(x)$ a $1$-parameter family of metrics on $\pl\bar{X}$ such that $h(0)=h_0$, 
with the regularity of $\rho v+\rho^n w$. 
We denote by $\bar{g_0}=x^2g_0$ and as before $\tra_{h_0}(h_1)=0$ if $h(x)=h_0+xh_1+O(x^2)$ (\emph{i.e.} $h_1$ is the second fundamental form of $\bar{g}_0$).
Then we can consider the elliptic Dirichlet problem (\ref{yamabe})
where $f\in C^{\infty}(\bar{X})$ is fixed.
It has a unique solution $u\in C^{\infty}(\bar{X})$ which
allows to define $\mc{N}:C^{\infty}(\pl\bar{X})\to C^{\infty}(\pl\bar{X})$ by
\[\mc{N}f=-\pl_n u|_{\pl\bar{X}}\]
where $u$ is the solution of \eqref{yamabe} and $\pl_n=\pl_\rho=\pl_x$ the interior unit normal 
vector field to $\pl\bar{X}$.
An easy computation as above shows that \eqref{yamabe} is equivalent to solving
\[\Big(\Delta_{g_0}-\frac{n+1}{2}.\frac{n-1}{2}\Big)\hat{\rho}^{\frac{n-1}{2}}u=0, \quad
u\in C^{\infty}(\bar{X}), \quad u|_{\pl\bar{X}}=f. 
\] 
by conformal covariance of the Yamabe operator.
But this is exactly the Poisson problem at energy $\la=(n+1)/2$ 
for the asymptotically hyperbolic manifold $g_0$, dealt with\footnote{They actually study it for smooth asymptotically hyperbolic manifolds but their proof works as well when 
log-terms enter the expansion of the metric at the boundary, in particular here the first $\log$ terms 
appear at order $x^{n}\log(x)$, thus they do not change the form of the two first asymptotic 
terms in the solution of the Poisson problem at energy $(n+1)/2$: one has 
\[\hat{\rho}^{\frac{n-1}{2}}u\sim x^{\frac{n-1}{2}}(f+O(x^2))+x^{\frac{n+1}{2}}\Big(-\mc{S}(\frac{n+1}{2})f+O(x)\Big)\]
where $\mc{S}(\la)$ is the scattering operator (see \cite{GZ}), the fact that there is no $x^{\frac{n+1}{2}}\log x$ terms
and no other terms than $x^{\frac{n+1}{2}}\mc{S}((n+1)/2)f$ at order $x^{\frac{n+1}{2}}$ 
is because $\tra_{h_0}(h_1)=0$ (then $\mc{S}(\la)$ has no residue at $(n+1)/2$), see Lemma 4.1
of \cite{G1} for more details.} by Graham-Zworski \cite{GZ}. 
Since $\pl_n=\pl_x$ we thus deduce that 
\[\mc{N}f=-S\Big(\frac{n+1}{2}\Big)f +\frac{(n-1)\omega}{2} f\]
where $\hat{\rho}=x(1+x\omega+O(x^2))$ for some $\omega\in C^{\infty}(\pl\bar{X})$.
But $\bar{g}x^2=\hat{\rho}^{2}\bar{g}_0$ thus applying to any vector field $V\in \pl\bar{X}$, this gives
\[(h_0(V)+xh_1(V))(1-2x\omega)=(h_0(V)+x \bar{h}_1(V))+O(x^2)\]
where $\bar{h}_1=\pl_{\rho}\bar{h}(0)$ is the second fundamental form of $\bar{g}$, this 
implies clearly that $\bar{h}_1=h_1-2\omega h_0$ and, taking the trace with respect to metric $h_0$,  
we get $\omega=-\frac{1}{2n}\tra_{h_0}(\bar{h}_1)$.
We have thus proved 
\begin{prop}\label{dnvsscat}
The Dirichlet-to-Neumann map for the conformal Laplacian of $\bar{g}$ is 
\[\mc{N}=\mc{S}\Big(\frac{n+1}{2}\Big)-\frac{(n-1)}{2}K\] 
where $K=\tra_{h_0}(\bar{h}_1)/2n$ is the mean curvature of $\pl\bar{X}$ for $\bar{g}$, $\mc{S}(\la)$
is the scattering operator associated to the complete metric with constant scalar curvature metric $g_0$, 
conformal to $\bar{g}$, and for choice of conformal representative $h_0=\bar{g}|_{T\pl\bar{X}}$.
\end{prop}

Another consequence, if $(\bar{X},\bar{g})$ is conformal to a convex
co-compact quotient $\Gamma\backslash\hh^{n+1}$ in even dimension, the
determinant of $\mc{N}+(n-1)K/2$ can be obtained by the functional equation
of \cite{G} in terms of a special values of Selberg zeta function of
$\Gamma$. We do not write the details and refer the reader to that paper,
since this is much less interesting than for surfaces.

\section{Appendix: the cylinder}

As a particular case of \cite{MT}, a smooth surface with boundary $\bar{X}$, 
with Euler characteristic $\chi(\bar{X})=0$, is conformal
to a hyperbolic cylinder $\cjg\gamma\cjd\backslash\hh^2$ with $\gamma:z\to e^{\ell} z$, which 
himself is conformal to the flat annulus $A_\rho:=\{z\in\cc;1<|z|<\rho\}$ with $\rho:=e^{2\pi^2/\ell}$,  
a conformal diffeomorphism being induced by the map
\[U:    z\in\mathbb H^2 \to  e^{2 i \pi (\log z)/\ell+2\pi^2/\ell}.\]
satisfying $U(z)=U(e^{\ell}z)$. 

We compute in this appendix the determinant of the Dirichlet-to-Neumann map for the annulus
$\overline{A_\rho}$ and check that if fits with the value found by the
technic of the functional equation of Selberg zeta function used above, giving
an alternative way of computing $\det'(\mc{N})$ for the cylinder, much in
the spirit of \cite{EW}.

Polar coordinates $z=re^{\idroit \theta}$ induce 
Fourier mode decompositions 
\[f=(f_\rho,f_1)\in 
L^2(\partial A_\rho)\simeq L^2(\rr/2\pi\zz)\otimes \mathbb C^2
\to (\widehat f_\rho(n),\widehat f_1(n))_{n\in\mathbb Z}\in \bigoplus_{n\in\mathbb Z}F_n
\]
with $F_n\simeq\mathbb C^2$ and the DN map  $\mc{N}\simeq
\bigoplus_{n\in\mathbb Z}\mc{N}_n$ is diagonal with respect to this
Fourier decomposition.
The harmonic functions $h_0$ and $h_1$ defined on $A_\rho$ by
\[
h_0(z)=1,\quad h_1(z)=1 -(1+\rho)\ln|z|/(\rho\ln\rho),\quad z\in A_\rho
\]
give eigenvectors of the DN map on $F_0$: $(1,1)$ and $(-\rho^{-1},1)$  with respective eigenvalues $0$
and $(1+\rho)/(\rho\ln\rho)$.

For $n\in\mathbb Z^*$, the harmonic functions
\begin{align*}
h_{\rho,n}(z)=\frac{z^n-\overline{z}^{-n}}{\rho^n-\rho^{-n}},
\quad 
h_{1,n}(z)=\frac{\rho^n\overline{z}^{-n}-\rho^{-n}{z}^{n}}{\rho^n-\rho^{-n}}
\end{align*}
have traces on $\partial A_\rho$ inducing the canonical base of
$F_n\simeq\cc^2$
\[
h_{\rho,n}(z)=\begin{cases}
  \eintheta&\hbox{if $|z|=\rho$}\\
  0&\hbox{if $|z|=1$}\\
\end{cases},\quad
h_{1,n}(z)=\begin{cases}
  0&\hbox{if $|z|=\rho$}\\
  \eintheta&\hbox{if $|z|=1$}\\
\end{cases}
\]
and their derivatives under the radial vector field $\partial_r$
\begin{align*}
\partial_rh_{\rho,n}(r\eintheta)=n\frac{r^{n-1}+r^{-n-1}}{\rho^n-\rho^{-n}}\eintheta
,\quad
\partial_rh_{1,n}(r\eintheta)=-n\frac{\rho^nr^{-n-1}+\rho^{-n}r^{n-1}}{\rho^n-\rho^{-n}}\eintheta
\end{align*}
give the matrix $\mc{N}_n$ with respect to the canonical base of $F_n$.
Observing that the interior normal derivative $\partial_n$ is $-\partial_r$
on $\{|z|=\rho\}$ and its opposite $\partial_r$ on $\{|z|=1\}$, we have,
with $\alpha=\log\rho$,
\begin{align}
\mc{N}_n=\begin{pmatrix}
\displaystyle n\rho^{-1}\frac{\rho^n+\rho^{-n}}{\rho^n-\rho^{-n}}&
\displaystyle -\frac{2n\rho^{-1}}{\rho^n-\rho^{-n}}\\
\displaystyle -\frac{2n}{\rho^n-\rho^{-n}}
&\displaystyle n\frac{\rho^n+\rho^{-n}}{\rho^n-\rho^{-n}}
\end{pmatrix}
=\frac{n}{\sinh(n\alpha)}
\begin{pmatrix}
\e{-\alpha}\cosh(n\alpha)&-\e{-\alpha}\\
-1&\cosh(n\alpha)
\end{pmatrix}
\label{eq:Mat}
\end{align}
with determinant $\delta_n=n^2\e{-\alpha}$ and
eigenvalues
\[ \lambda_{n,\pm}=\frac{\e{-\alpha/2}n}{\sinh(n\alpha)}{\cosh(\alpha/2)\cosh(n\alpha)\pm\sqrt{\sinh^2(\alpha/2)\cosh^2(n\alpha      )+1}}.
\]

The following lemma claims the product relation
``${\prod \lambda_{n,+}\prod \lambda_{n,-}=\prod \delta_n}$'':
\begin{lem}
Let $u=(u_n)_{n\ge1}$ a sequence  with positive terms such that
$u_n=n^k(1+\varepsilon_n)$ with $\varepsilon_n={O}(\e{-a n})$ for
some non negative $k$ and $a$. If $\zeta_u$ is the zeta function
$\zeta_u(s)=\sum_{n\ge1}u_n^{-s}$ convergent for $\Re s >k^{-1}$, then $\zeta_u$
extends meromorphically to the complex plane and is regular for $s=0$ with
\[
\zeta_u(0)=\zeta(0),\quad
\partial_s\zeta_u(0)=k\zeta'(0)-\sum_{n\ge1}\ln(1+\varepsilon_n)
\]
where $\zeta$ is the Riemann zeta function. If $u,v $ are two such
sequences, then
the sequence $w$ defined by $w_n=u_nv_n,n\ge1$ is again of the same type
and $\partial_s\zeta_u(0)+\partial_s\zeta_v(0)=\partial_s\zeta_w(0)$.
\end{lem}
\Proof Let $G(s,\varepsilon)$ be the function defined for $\varepsilon$
small and $s\in\mathbb C$, holomorphic in $s$, such that
\[
(1+\varepsilon)^{-s}=1-sG(s,\varepsilon),
\quad |G(s,\varepsilon)|+ |\partial_sG(s,\varepsilon)|=_{s,\varepsilon\sim0}{O}(\varepsilon),
\quad G(s,\varepsilon)_{|s=0}=\ln(1+\varepsilon).
\]
We have then
\[
\zeta_u(s)=\sum_{n\ge1}(n^k(1+\varepsilon_n))^{-s}
=\sum_{n\ge1}n^{-ks}-s \sum_{n\ge1}n^{-ks}G(s,\varepsilon_n)
=\zeta(ks)-s \sum_{n\ge1}n^{-ks}G(s,\varepsilon_n)
\]
and
\[
\zeta_u(0)=\zeta(0),\quad
\partial_s\zeta_u(0)=k\zeta'(0)-\sum_{n\ge1}G(0,\varepsilon_n)
=k\zeta'(0)-\sum_{n\ge1}\ln(1+\varepsilon_n).\qedhere
\]
\medskip The $\zeta$-regularized determinant
$\det'\mc{N}$ for the DN map on the annulus
$\overline{A_\rho}$ is defined through the zeta function
\[
\zeta_\rho(s)=[(1+\rho)/(\rho\ln\rho)]^{-s}+\sum_{n\in\mathbb Z^*}
\left [\lambda_{-,n}^{-s}+\lambda_{+,n}^{-s}\right].
\]
According to the preceding Lemma,
if $\widetilde\zeta_\rho$ is the zeta function defined by  
$\widetilde\zeta_\rho(s)=[(1+\rho)/(\rho\ln\rho)]^{-s}+2\sum_{n\ge1}
\delta_n^{-s}$.  We have
\[\partial_s\zeta_\rho(0)= \partial_s\widetilde\zeta_\rho(0)
=-\ln[(1+\rho)/(\rho\ln\rho)]
+2\alpha\zeta(0)+4\zeta'(0)
=-\ln[(1+\rho)/(\ln\rho)] -2\ln(2\pi),\]
where we have used   $\la_{n,+}\la_{n,-}=\delta_{|n|}=n^2e^{-\alpha}$ (so
that the $\log(1+\eps_n)$ terms disappear in the Lemma) 
and for the last equality $\zeta(0)=-1/2$ and $\zeta'(0)=-\ln(2\pi)/2$. Hence
\[
 \frac{\hbox{$\det'$}\mc{N}}{\ell(\partial
\overline{  A_\rho})}
=\frac{\e{-\partial_s\zeta_\rho(0)}}{2\pi(1+\rho)}
=\frac{2\pi}{\ln\rho}=\frac{\ell}{\pi}
\]
\bigskip
which perfectly fits with Theorem \ref{detDN}.

\end{document}